\documentclass[11pt,a4paper]{amsart}

\usepackage{latexsym}

\newtheorem{theorem}{Theorem}[section]
\newtheorem{lemma}[theorem]{Lemma}
\newtheorem{prop}[theorem]{Proposition}

\theoremstyle{definition}

\theoremstyle{remark}

\renewcommand{\theequation}{\arabic{section}.\arabic{equation}}

\def\Q{{\mathbb{Q}}}  \def\C{{\mathbb{C}}}
\def\Z{{\mathbb{Z}}}  \def\P{{\mathbb{P}}}
 
 \def\Qaf{\overline{\mathbb{Q}}}
\def\Cal{\mathcal}

\newcommand{\Gal}{\operatorname{Gal}}

\newcommand{\rank}{\operatorname{rank}}

\newcommand{\End}{\operatorname{End}}

  \def\om{\omega} 

\def\phihat{\kern.3em\hat{\kern-.3em\phi}}
\def\sphihat{\kern.2em\hat{\kern-.2em\phi}}

\def\<{\langle} \def\>{\rangle}

\makeatletter
\def\currenteq{\textup{\tagform@{\p@equation\theequation} }}
\makeatother

\title
{Sections on certain $j=0$ elliptic surfaces}

\author{Jasbir Chahal}
\address{Department of Mathematics \\
Brigham Young University \\
Provo, UT 84602-6539\\
USA}
\email{jasbir@math.byu.edu}

\author{Matthijs Meijer}
\author{Jaap Top}
\address{Vakgroep Wiskunde R{\sl u}G \\
P.O. Box 800 \\
9700 AV Groningen \\
the Netherlands}
\email{top@math.rug.nl}

\thanks{This work originated from a lecture by Masato Kuwata
during the 914th AMS meeting at Rider university, Lawrenceville,
New Jersey in October 1996. The third author would like to thank
Henri Darmon and Fernando Rodriguez Villegas for inviting him
to this meeting, and Masato Kuwata and Jasper Scholten for
their interest in this work. Several ideas in this paper were
developed during the summer of 1998 which the third author
spent at BYU. He thanks Jasbir Chahal, Bill Lang and Tyler Jarvis
for their hospitality. The middle section of this paper was a
part of the second author's Master's thesis in mathematics at the
university of Groningen.}

\date{November 1999}
\begin{document}

\maketitle

\section{Introduction}
It is well known that elliptic curves $E$ over a function field $k(t)$
are in 1-1 correspondence with elliptic surfaces
${\Cal E}\rightarrow\P^1$ over $k$. The geometric interpretation
of points in $E(k(t))$ arising in this way has led to a lot
of results; see, e.g., \cite{Shi72}, \cite{Shi90}.
For instance, if $\mbox{char}(k)=0$ and $E$ is not isomorphic 
over $k(t)$
to a curve already defined over $k$, then
$\rank\,E(k(t))\leq 10p_g+8$, where $p_g$ is the geometric
genus (the number of independent regular 2-forms) of the elliptic 
surface.

In case the surface is rational one has $p_g=0$ and hence the
upper bound above equals $8$. The pencil of cubic curves passing
through $8$ given points (in general position) in the plane provides 
an example of such a surface. The $9$ base points of the pencil
yield sections of the surface, and with any one of them as
origin the others are (in the general case) independent. More details 
on this are found in \cite{Shi90} and references given there. In fact
it is not hard to construct rational surfaces
with any rank between $0$ and $8$. In particular explicit examples
exist even with a basis for the group $E(k(t))$.

It is interesting to note that the upper bound $10p_g+8$
is not known to be sharp for high geometric genus. In fact
it is not even known whether $\rank\,E(\C(t))$ can attain
arbitrarily high values. The current record seems to be due
to Shioda (1992). He asserts in \cite{Shi92} that 
the curve given by $y^2=x^3+t^{360}+1$ has rank $68$ over
$\C(t)$. This example has $p_g=59$ hence the actual rank
is much smaller than the general theoretical upper bound.
We will discuss this example in more detail below (proposition~\ref{rk36}).

The smallest case beyond rational surfaces is the one with
$p_g=1$, so ranks $\leq 10\cdot 1+8=18$. The elliptic surfaces 
in this case are so-called $K3$-surfaces. A result of David Cox states
that any integer $r$ with $0\leq r\leq 18$ occurs as the
rank of some elliptic $K3$-surface over $\C$; see \cite{Cox} and also
\cite{Nish}. The proof uses the surjectivity of a period map,
and is therefore transcendental in nature. In particular,
it does not give equations of such surfaces, let alone
explicit descriptions of independent sections. A different
proof for the existence of elliptic $K3$-surfaces with 
Mordell-Weil rank $18$ is presented in Kuwata's paper \cite{Kuw}.
It is based on work of Inose which implies, that if a rational
map of finite degree exists between $K3$-surfaces, then
they have the same N\'{e}ron-Severi rank (\cite[Cor.~1.2]{In}). 
Kuwata starts
with the Kummer surface of a product $E_1\times E_2$ of
isogenous CM-elliptic curves. He then constructs a degree $3$
cover of it which is still $K3$, and which moreover in many cases 
admits an elliptic fibration without reducible fibres.
This easily implies that the resulting example has Mordell-Weil
rank $18$. It is straightforward to write down explicit
examples of such surfaces, even with equations over $\Q$ (although
the sections will usually only be defined over an extension).
However, from the transcendental nature of the result of Inose,
it remains unclear how to find independent sections for
these examples. 

In the present paper we present an example with an explicit
description of $18$ independent sections:
\begin{theorem}\label{main}
The elliptic curve $E/\Q(t)$, defined by the equation
$y^2=x^3-27(t^{12}-11t^6-1)$ corresponds to a $K3$-surface.
It has $\rank\,E(\Qaf(t))=18$.

One finds $18$ independent
points among the $(x(t),y(t))\in E(\Qaf(t))$ where $x(t)$
has one of the forms 
$$\frac{\displaystyle at^{12}+bt^6+c}{\displaystyle t^n}$$
 or 
$$\frac{\displaystyle \om \left(
t^{24}-12t^{18}+14t^{12}+12t^6+1\right)}{\displaystyle 4t^{10}},$$
with $n\in\{0,2,4,6,8\}$ and $\om^3=1$.
 
\end{theorem}

It seems plausible that one can prove this by merely finding
sufficiently many points as in the theorem, and then
calculating a height pairing determinant. However, we
have not tried this rather uninteresting and elaborous method.
Instead, we give a proof which is hopefully more
enlightening. It relies on the fact that the elliptic
curve given in the theorem has $j$-invariant $0$. Some Galois
theory and linear algebra then allow one to break the
Mordell-Weil group into parts. Most of these parts
correspond to Mordell-Weil groups of rational elliptic
surfaces, which are well understood. 

We finish the paper by applying the same method to the
problem of constructing high rank curves over function fields.
In particular, the
rank $68$ example of Shioda is studied, and we mention some examples
of elliptic curves of a similar form with a fairly high
$\Q(t)$-rank.

\section{Generalities on $j=0$ elliptic curves}
Suppose $E$ is an elliptic curve defined over a field $K$.
For any finite extension $L/K$, the Galois group $\Gal(L/K)$
acts on the group $E(L)$ of $L$-rational points of $E$.
Regarding $E(L)$ as a module over $R:=\End_K(E)$, this
Galois action is $R$-linear. This can be used to
decompose $E(L)$, or rather the vector space $E(L)\otimes_{\Z}\Q$,
as a direct sum of smaller vector spaces.

The standard example where this is used, is given by
a quadratic Galois extension $L/K$. The Galois action splits
$E(L)\otimes_{\Z}\Q$ as $V^+\oplus V^-$. Here $V^+$ consists
of the elements in $E(L)\otimes_{\Z}\Q$ on which $\Gal(L/K)$
acts trivially; hence $V^+=E(K)\otimes_{\Z}\Q$. The other
space $V^-$ consists of the points in 
$E(L)\otimes_{\Z}\Q$ on which $\Gal(L/K)$ acts via the
non-trivial character of order 2. In other words,
if $\tau\in \Gal(L/K)$ is the generator, then $V^-$
consists of all $P\otimes r$ for which $\tau(P)=-P$.
This implies that $V^-$ can be identified with
$E'(K)\otimes_{\Z}\Q$, in which $E'$ is the quadratic
twist of $E/K$ over $L$. In particular, if we are in
a situation where the Mordell-Weil theorem holds, then
$$\rank\,E(L)=\rank\,E(K)+\rank\,E'(K)$$
(compare \cite[exerc. X-10.16]{Sil}) as is well known.

We will use such a splitting in the case that $\End_K(E)$
has field of fractions $\Q(\om)$, for a primitive cube root
of unity $\om$. Suppose that we are in this case.
Then for any extension $L/K$, the vector space 
$E(L)\otimes_{\Z}\Q$ is a linear space over $\Q(\om)$.
In particular, if $L/K$ is a Galois extension with a
cyclic Galois group $\langle \sigma\rangle$ of order $6$,
then $\sigma$ acts $\Q(\om)$-linearly on $E(L)\otimes_{\Z}\Q$.
Since all eigenvalues $\lambda$ of $\sigma$ satisfy
$\lambda^6=1$, they are $\Q(\om)$-rational, hence
$E(L)\otimes_{\Z}\Q$ splits as a direct sum $\sum V_{\lambda}$.
Here $V_{\lambda}$ is the $\Q(\om)$-vector space consisting of 
the points $P\otimes r\in E(L)\otimes_{\Z}\Q$ for which 
$\tau(P)=\lambda P$.
As in the standard example, each $V_\lambda$ can be
interpreted in terms of $K$-rational points on a certain
(quadratic or cubic or sextic) twist of $E/K$ over $L$.

This is made explicit in the following example.
\begin{lemma}\label{ranksum}
Let $k/\Q(\om)$ be a field extension and let $K=k(s)$ be a purely
transcendental extension of $k$. Consider 
$L:=k(t)$ over $K$, in which $t^6=s$. Then $L/K$ is cyclic
of degree $6$.

If $E/K$ is the elliptic curve defined by the equation
$y^2=x^3+f(s)$ for some rational function $f(s)\in K$,
then the $\Q(\om)$-vector space $E(L)\otimes_{\Z} \Q$
splits as a direct sum 
$$\sum_{n=0}^5 V_{(-\om^2)^n}.$$
The $\Q(\om)$-vector space $V_{(-\om^2)^n}$ here can
be identified with $E_n(K)\otimes_{\Z}\Q$, in which $E_n/K$
is the elliptic curve given by $y^2=x^3+s^nf(s)$. 
\end{lemma}

\noindent
{\it Proof.}
Fix $\End_K(E)\cong \Z[\om]$ by
identifying $\om$ with the endomorphism $(x,y)\mapsto(\om x,y)$.
The only statement of the lemma which is not explained in
the discussion above, is the identification of the
eigenspace $V_{(-\om^2)^n}$ with $E_n(K)\otimes_{\Z}\Q$.
This is a direct consequence of the fact that 
$P=(x(t),y(t))$ is in $E(L)$ and satisfies $(x(-\om^2t),y(-\om^2t))
=(-\om^2)^n P$, precisely when $(t^{2n}x(t),t^{3n}y(t))\in E_n(K)$.
This shows in particular how a point $P\otimes r\in
V_{(-\om^2)^n}$ yields an element of $E_n(K)$. Conversely, 
$(x(s),y(s)\in E_n(K)$ is sent to 
$\left(x(t^6)/t^{2n},y(t^6)/t^{3n}\right)\in E(L)$ in the associated 
eigenspace. \hfill{$\Box$}

\par
In the rest of this paper we only consider $j=0$ elliptic curves $E/K$,
given by an equation $y^2=x^3+k$. Note that for $\om\in K$ the 
$\Z[\om]$-module
structure on $E(K)$ implies, that if $P\in E(K)$ is a point of
infinite order, then so is $\om P$, and $P,\om P\in E(K)$ are
independent (over $\Z$). More generally, the fact that the
$\Z$-module $E(K)$ is even a $\Z[\om]$-module, shows that if
it has finite rank over $\Z$, then this rank is even.
Some related information is given in the
following lemma.
\begin{lemma}\label{qrank}
Suppose that $\rank\,E(K)$ is finite, and that $L:=K(\sqrt{-3})$ is
a quadratic extension of $K$. Then $\rank\,E(L)=2\rank\,E(K)$.
\end{lemma} 

\noindent{\sl Proof.} One way to see this, is to split up $E(L)$
for the action of $\Gal(L/K)$. The twist $E'/K$ that is associated to
the $-1$-eigenspace is over $K$ isogenous (by an isogeny of degree $3$)
to $E$. Hence $E(K)$ and $E'(K)$ have the same rank, and the
sum of these two ranks is the rank of $E(L)$.

A somewhat different proof runs as follows. Let $\tau$ denote the
non-trivial element of $\Gal(L/K)$, and define $\phi$ on $E(L)$
by $\phi(P)=\om P$. Consider the sequence
$$0\rightarrow E(K)\stackrel{\phi}{\longrightarrow} E(L)
\stackrel{\phi^2+\tau\phi^2}{\longrightarrow} E(K).$$
It is easily verified that this is an exact sequence; e.g.,
use that $\phi^2+\tau\phi^2=\phi^2+\phi\tau$ on $E(L)$.
Moreover, the cokernel of the rightmost arrow is torsion,
because $(1+\tau)\phi^2\phi^4E(K)=2E(K)$. This implies the lemma.
\hfill{$\Box$}

\par Note that this lemma implies for the curve given in
theorem~\ref{main}, that although it has rank $18$ over $\Qaf(t)$,
its rank over $\Q(t)$ is at most $9$. In fact this $\Q(t)$-rank is
much smaller than $9$, as can be calculated using \cite{Br} or
alternatively, by computing a characteristic polynomial of
Frobenius acting on the second $\ell$-adic cohomology group
of the elliptic surface over $\Qaf$ corresponding to $E$;
see \cite{vG-T}, \cite[Ch.~2]{Sch} for details of such a
computation. 

\par
Combining the two lemmas, one concludes that for any 
extension $k(s)=k(t^6)\subset k(t)$ of degree $6$ and elliptic curve 
$E/k(s)$ given
by $y^2=x^3+f(s)$, one has
$$\rank\,E(k(t))=\sum_{n=0}^5\rank E_n(k(s)),$$
at least if these ranks are finite. As before, here $E_n$
denotes the elliptic curve given by $y^2=x^3+s^nf(s)$.
In particular this formula holds for
all number fields $k$, whether or not $k$ contains a primitive
cube root of unity. 
\section{The rank $18$ example}
We will now present a proof of theorem~\ref{main}. In fact, it
will be explained at the same time how to construct this
example, and even that it is unique in a certain sense.

Consider a polynomial $f(s)\in\Qaf[s]$ of degree $d>0$, which
is not divisible by a $6$th power. We assume that $f(0)\neq 0$;
this implies that the polynomial $f(t^6)\in\Qaf[t]$ is
$6$th power free as well. The elliptic surface over $\Qaf$
corresponding to the equation $y^2=x^3+f(t^6)$ then has
geometric genus $p_g=d-1$. In particular, one obtains a $K3$-surface
precisely when $f$ has degree $2$. Hence from now on
$E$, given by the equation $y^2=x^3+at^{12}+bt^6+c$, with
$a,b,c\in\Qaf$ satisfying $ac\neq 0$, is considered.
The elliptic fibration associated with this Weierstrass equation has
singular fibres exactly over the roots of $ at^{12}+bt^6+c=0$.
If $b^2-4ac=0$, then this has $6$ zeroes, each with
multiplicity $2$. This implies that in that case the elliptic
$K3$-surface has $6$ singular fibres, each of Kodaira type $IV$.
The Shioda-Tate formula then yields the upper bound
$18-6\cdot 2=6$ for $\rank E(\Qaf(t))$. We will assume
that this situation does not occur; in other words, $b^2-4ac\neq 0$.

By lemma~\ref{ranksum}, $\rank E(\Qaf(t))=\sum_{n=0}^5\rank
E_n(\Qaf(s))$. In the present case, $E_n$ is given by
$y^2=x^3+s^n(as^2+bs+c)$. For $0\leq n\leq 4$ the corresponding
elliptic surface is a rational surface. For such surfaces
the Shioda-Tate formula provides an exact formula for
$\rank E_n(\Qaf(s))$ and it is even known that
the group is generated by points whose $x$-coordinate
is a polynomial of degree at most $2$. See \cite{Shi90} for
details. For $n=0, n=4$ the only reducible fibre of the
associated elliptic surface is of type $IV^*$, which implies
that the Mordell-Weil rank equals $2$. For $n=1, n=3$ one
again finds only one reducible fibre; it is of type $I_0^*$
hence here the rank is $4$. In the remaining case $n=2$ there
are $2$ reducible fibres, each of type $IV$. This gives
$4$ as the rank of $E_2$ over $\Qaf(s)$.
Adding these contributions and using the proof of
lemma~\ref{ranksum} one now has $16$ independent points in
$E(\Qaf(t))$, all of the form as asserted in the statement
of theorem~\ref{main}.

It remains to consider the only eigenspace not taken into
account yet, namely the one corresponding to $E_5$. The
following well known lemma will be used.
\begin{lemma}\label{notors}
Let $k$ be a field of characteristic $0$, and let $E/k(s)$
be an elliptic curve given by $y^2=x^3+g(s)$. Suppose that
$g(s)=0$ has some simple root in an algebraic closure of $k$.
Then $E(k(s))$ is a finitely generated free group. 
\end{lemma} 

\noindent{\sl Proof.} Since $g(s)$ has a simple zero, $E$ is
not isomorphic over $k(s)$ to an elliptic curve already defined
over $k$. This implies that the Mordell-Weil theorem holds
for $E(k(s))$; in other words, it is a finitely generated
group. To show that it is torsion free, we specialize at
the simple root. This defines a map $E(\overline{k}(s))\rightarrow
C(\overline{k})$ where $\overline{k}$ denotes an algebraic
closure of $k$ and $C$ denotes the curve given by $y^2z=x^3$ in $\P^2$.
Since we specialize at a simple root, the singular point on $C$
is not in the image, and the specialization is in fact a homomorphism
of groups. The target group is the additive group $\overline{k}$.
Hence since such a specialization is known to be injective
on torsion points, the lemma follows. \hfill{$\Box$}

\par
Returning to the proof of theorem~\ref{main}, all we have to do
in order to have an example of rank $18$ is to exhibit a non-trivial
point on $E_5$, which is the curve given by $y^2=x^3+as^7+bs^6+cs^5$.
We will try to find such a point of the simplest possible kind,
namely with $x$-coordinate a polynomial in $\Qaf[s]$. It is clear
from the equation that such a polynomial should have degree at 
least $4$. We therefore look at the smallest possible degree,
namely $4$. The corresponding $y$-coordinate should then be a
polynomial of degree $6$. We can and will assume that the
leading coefficients of these coordinates are equal, and then by
scaling $s$ if necessary, that they equal $1$.
This leads to the problem of finding $a_0,\ldots,a_3$
and $a,b,c$ and $b_1,\ldots,b_5$ in $\Qaf$ such that
$$(a_0^3+b_1s+b_2s^2+\ldots+b_5s^5+s^6)^2=(a_0^2+a_1s+a_2s^2+a_3s^3+s^4)^3
+as^7+bs^6+cs^5.$$

Using Maple, and eliminating as many as possible of the variables which
appear linearly when comparing coefficients, this problem can be
solved completely. Apart from solutions with $a=b=c=0$, and multiplying
$a,b,c$ with some given sixth power (which does not change the
elliptic curve $E$), one finds the unique additional solution
$(a,b,c)=-27\cdot(1,-11,-1)$. The corresponding degree $4$ polynomial
as $x$-coordinate is upto a third root of unity
$x(s)=(s^4-12s^3+14 s^2+12 s+1)/4$.
This implies the theorem. Note that in fact
$\rank E_5(\Qaf(s))=2$ since it is $\geq 1$ and even by what we showed
above and by the $\Z[\om]$-structure, and it is $\leq 2$ because
of the formula $\sum\rank E_n(\Qaf(s))\leq 18$. \hfill{$\Box$}

\vspace{\baselineskip}
\noindent{\sl Remark.}
It may be of some interest to relate the above construction to
the results of Inose and of Kuwata which were mentioned in the
introduction.

Start with the $E_5$ above that has $\Qaf(s)$-rank
at least $2$. This curve corresponds to a $K3$-surface over
$\Qaf$, and the given elliptic fibration on it has precisely
$2$ reducible fibres, both of type $II^*$ (at $s=0$ and at 
$s=\infty$). The Shioda-Tate formula for the rank $\rho$ of
the N\'{e}ron-Severi group of this surface now implies
$$20\geq \rho \geq 2+2\cdot 8+\rank E_5(\Qaf(s)) \geq 20,$$
hence we conclude (as before) that 
$\rank E_5(\Qaf(s))=2$ and also that $\rho=20$.

The result of Inose \cite{In} states, that any $K3$-surface
which admits a finite rational map to our one, will then
also have $\rho=20$. An example of such a surface is the
one corresponding to $E$; a rational map in this case is 
given by
$(x,y,t)\mapsto (x/t^{10},y/t^{15},s=t^6)$. Since the
latter surface has no reducible fibres, the Shioda-Tate
formula implies that $\rank E(\Qaf(t))=18$ as we
already proved.

As a special case of a construction of Kuwata \cite{Kuw},
take the elliptic curves with equations
$y^2=x^3+1$ and $y^2=x^3-15x+22$. These curves are $2$-isogenous,
complex multiplication curves. This implies that the Kummer surface
of their product is a $K3$-surface with a N\'{e}ron-Severi group
of rank $20$. An affine equation for this Kummer surface is
$(t^3+1)y^2=x^3-15x+22$. Now the surface $X$ corresponding to
$(t^3+1)\eta^6 =x^3-15x+22$ admits a rational map of degree $3$
to this Kummer surface (given as $\eta\mapsto y=\eta^3$).
Note that $X$ admits an elliptic fibration $(x,\eta,t)\mapsto \eta$,
with for instance $\eta\mapsto (x=2,\eta,t=-1)$ as a section.
One computes that a Weierstrass equation for this fibration
is $Y^2=X^3+27\cdot16((2\eta^3)^4+11(2\eta^3)^2-1)$. Hence the
surface $X$ is over $\Qaf$ isomorphic to the one presented
in theorem~\ref{main}. Note that Inose's result implies that
it has N\'eron-Severi rank $20$, hence a third proof for
the fact that $\rank E(\Qaf(t))=18$ is obtained. This
method of Kuwata has the advantage that it works for other pairs
of isogenous CM-curves as well, but only in the present case
 $18$ independent points have been found.

\section{Other examples}
Two aspects of the methods above will be discussed: how to apply
them in order to find even higher ranks over $\Qaf(t)$, and what
can be found over $\Q(t)$.

\subsection{Higher ranks}
If we are not interested in $K3$-surfaces only, it is possible
to give examples of the form $y^2=x^3+f(t^6)$ with a higher rank.
The following results illustrate this.
\begin{prop}
Suppose $f(s)=s^3+as^2+as+1\in\Qaf[s]$ is a polynomial with simple
zeroes. Then $E/\Qaf(t)$, given by $y^2=x^3+f(t^6)$
satisfies $\rank E(\Qaf(t))\in\{20,24,28\}$. Moreover, there exists
$a\in\Qaf$ for which $\rank E(\Qaf(t))\geq 24$.  
\end{prop}

\noindent{\sl Proof.} We use lemma~\ref{ranksum} and the notations
introduced there. For a cubic $f(s)$ as above, the curves
$E_0,E_1,E_2$ and $E_3$ correspond to rational elliptic surfaces.
Adding their contributions yields a subgroup of rank $20$
in $E(\Qaf(t)$. For the remaining two eigenspaces, note that
$(x,y,s)\mapsto (x/s^4,y/s^6,1/s)$ defines an isomorphism between
$E_4$ and $E_5$. The elliptic $K3$-surface corresponding to $E_4$
has two reducible fibres of type $IV^*, II^*$ respectively.
Hence $\rank E_4(\Qaf(s))\leq 18-8-6=4$. Moreover, we have seen
that this rank is even. Since
$\rank E(\Qaf(t))=20+2\rank E_4(\Qaf(s))$ this implies the three
mentioned possibilities for the rank.

To show that rank $\geq 24$ occurs, it suffices by the above
and lemma~\ref{notors} to find a value $a\in\Qaf$ for which
$E_4(\Qaf(s))\neq \{O\}$. This can be done as in the proof
of theorem~\ref{main}. \hfill{$\Box$}

\par
In fact much higher ranks can be obtained by iterating the method
of decomposing the space of sections into eigenspaces.
Denote by $r(f(t))$ the rank of the group $E(\Qaf(t))$, where
$E$ is given by $y^2=x^3+f(t)$. Using this notation,
lemma~\ref{ranksum} gives the formula
$$r(f(t^6))= r(f(t))+r(tf(t))+r(t^2f(t))+r(t^3f(t))+r(t^4f(t))+r(t^5f(t)).$$
If we only use the field automorphism of $\Qaf(t)$ that multiplies
$t$ by a third or by a second root of unity, the analogous
formulas
$$r(f(t^3))=r(f(t))+r(t^2f(t))+r(t^4f(t))$$
and
$$r(f(t^2))=r(f(t))+r(t^3f(t))$$
are obtained. Starting with any polynomial $g(t)$, the above three
formulas can be applied repeatedly until none of the resulting 
polynomials is a polynomial in $t^2$ or in $t^3$. As an example,
one could start with a polynomial $g$ with $g(0)\neq 0$
and take $f(t)=g(t^n)$. Decomposing as far as possible
leads to a certain set of polynomials $t^ag(t^b)$, with $0\leq a\leq 5$,
and $\gcd(a,b)=1$ and $b$ is a divisor of $n$ such that $n/b$ is a power
of $2$ times a power of $3$. This results in a lower bound for
$r(f(t))$ if we take only those contributions $r(t^ag(t^b))$ that
correspond to rational elliptic surfaces. In other words,
we only consider the $a,b$ for which $a+b\cdot\deg(g)\leq 6$.

The maximal number of such possible contributions is obtained
when one considers $\deg(g)=1$. In this case, in order to have all
terms $r(t^ag(t))$, we need that $n$ is a power of $2$ times
a power of $3$ and moreover $6|n$. To have the terms
$r(tg(t^2))$ and $r(t^3g(t^2))$ as well leads to the additional
assumption $12|n$. Next, the terms $r(tg(t^3))$ and $r(t^2g(t^3))$
are present if also $18|n$. And finally the term $r(tg(t^4))$
appears in addition by demanding $24|n$.  The result one
obtains in this way, is the following.
\begin{prop}\label{rk36}
Let $n$ be a positive multiple of $72$ and let $g(t)\in\Qaf[t]$ be a degree $1$
polynomial with $g(0)\neq 0$.

The curve given by $y^2=x^3+g(t^n)$ has rank $\geq 36$ over $\Qaf(t)$.
Moreover, one can find $36$ independent points all coming from rational
elliptic surfaces via a base change.

If $n$ is  a positive multiple of $360$, then the rank of this curve
over $\Qaf(t)$ equals $68$. Here $60$ independent points come  from rational
elliptic surfaces via a base change, and $8$ more in the same way from
an elliptic $K3$-surface.
\end{prop}   

\noindent{\sl Proof.} As already remarked in the
introduction of this paper, the rank $68$ example here is
due to Shioda. Using his algorithm \cite[Thm.~1]{Shi86},
one can in particular calculate the rank of any curve as given
in the proposition; the result is that the
rank is $\leq 68$, with equality if and only if $n$ is a positive
multiple of $360$.
We will now follow the discussion
above. This will lead to $36$ resp. $60$ independent points coming from
the various associated rational elliptic surfaces. Note that
in particular, such points may be explicitly given.

To prove the first part of the proposition, note that $r(g(t^n))\geq r(g(t^m))$
if $m|n$. Hence we may and will assume $n=72$. Then $r(g(t^{72}))$\\
$
\begin{array}{rl}
\mbox{  }\hspace{12pt}&\geq  r(g(t^6))+r(tg(t^2))+r(t^3g(t^2))+r(tg(t^3))+
r(t^2g(t^3))+r(tg(t^4))\\
& = 8+4+4+6+6+8=36,
\end{array}
$
which is what we wanted to prove. 
\par
It remains to prove the assertions about the case of
a positive 
$n\equiv 0\bmod 360$.
In the discussion preceding proposition~\ref{rk36}, we
only considered rational eigenspaces corresponding to
$g(t^6)$, $tg(t^4)$, $tg(t^3),t^2g(t^3)$ and
$tg(t^2),t^3g(t^2)$. One can use the
remaining ones, for $g(t^5)$ and $tg(t^5)$ as well,
by relaxing the assumption that some eigenspace corresponds to
$y^2=x^3+t^ag(t^b)$ to the weaker statement that the points on such
a curve give rise to a subspace (of possibly smaller
dimension) in an eigenspace. We illustrate this by
an example. Starting from the polynomial $g(t^{30})$,
one of the eigenspaces when using the automorphism of order
$6$ will correspond to $t^5g(t^5)$. Using our automorphisms
of order $2$, $3$ or $6$, we cannot decompose this space
into smaller pieces. However, the invariants for the order $5$
automorphism that multiplies $t$ by a fifth root of unity,
will give a subspace corresponding to $sg(s)$. Hence
$r(t^5g(t^5))\geq r(tg(t))$, and the latter number may
correspond to a rational surface even though the first one
does not. 
Although it will not be used here, we remark that such a
difference $r(t^5g(t^5))- r(tg(t))$ is a multiple of $8$.
To prove this, note that we consider the $\Q(\om)$-vector space
$V=E(\Qaf(t))\otimes_{\Z}\Q$ here, with $E$ given by
$y^2=x^3+t^5g(t^5)$. This vector space comes equipped
with an action of a cyclic group of order $5$. Since
a nontrivial fifth root of unity generates an extension
of degree $4$ over $\Q(\om)$, the representation on $V$ splits
into a sum $V_{\mbox{\scriptsize tr}}\oplus V_{\mbox{\scriptsize nt}}$,
where the action on $V_{\mbox{\scriptsize tr}}$ is the trivial one
and where $V_{\mbox{\scriptsize nt}}$ is a sum of $4$-dimensional
irreducible representations. In particular,
$\dim_{\Q(\om)}V_{\mbox{\scriptsize nt}}$ is a multiple
of $4$, and therefore
$r(t^5g(t^5))- r(tg(t))=\dim_{\Q}V_{\mbox{\scriptsize nt}}$
is a multiple of $8$.

We now discuss a number of eigenspaces that come up
when starting with a polynomial $g(t^{360})$. Here $g$ is
a polynomial of degree $1$. Throughout,
the assumption $g(0)\neq 0$ is made. Only in the
spaces associated with $g(t^{60}), t^2g(t^{60}), t^3g(t^{60})$
and $t^4g(t^{60})$ we are able to find pieces
corresponding to rational elliptic surfaces.
\begin{enumerate}
\item[$t^2g(t^{60})$.] Using order $2$ automorphisms, this part
yields 
$$\begin{array}{rl}
r(tg(t^{30}))+r(t^4g(t^{30})) &\geq r(t^2g(t^{15}))+ r(t^5g(t^{15}))\\
& \geq r(tg(t^3))= 6.\end{array}
$$
\item[$t^4g(t^{60})$.] Analogous to the previous part, one finds here
$$\begin{array}{rl}
r(t^2g(t^{30}))+r(t^5g(t^{30})) &\geq r(tg(t^{15}))+ r(t^4g(t^{15}))\\
&\geq r(t^4g(t^{15}))=r(t^6\cdot t^4g(t^{15}))\\
&\geq r(t^2g(t^3))=6.\end{array}
$$ 
\item[$t^3g(t^{60})$.] Here a splitting into three parts can be used;
it gives
$$\begin{array}{rl}
r(tg(t^{20}))+r(t^3g(t^{20}))+r(t^5g(t^{20}))& \geq r(t^5g(t^{20}))\\
& \geq r(tg(t^4))=8.
\end{array}
$$ 
\item[$g(t^{60})$.] In this case one obtains a splitting into $6$ parts,
corresponding to the polynomials $t^ng(t^{10})$ for $0\leq n\leq 5$.
These parts will now be analyzed separately.
\begin{enumerate}
\item[$n=5$.] We estimate $r(t^5g(t^{10}))\geq r(tg(t^2))=4$.
\item[$n=4$.] Analogously, $r(t^4g(t^{10}))=r(t^{10}g(t^{10}))
\geq r(t^2g(t^2))=4$.
\item[$n=3$.] This gives $r(t^3g(t^{10}))=r(t^{15}g(t^{10}))
\geq r(t^3g(t^2))=4$.
\item[$n=2$.] Using a further splitting one obtains
$$\begin{array}{rl}
r(t^2g(t^{10}))&= r(tg(t^5))+r(t^4g(t^5))\\
&=r(tg(t^5))+r(t^{10}g(t^5))\\
&\geq r(tg(t^5))+r(t^2g(t))=8+2=10.
\end{array}
$$
\item[$n=0$.] Just as in the previous case,
$$\begin{array}{rl}
r(g(t^{10}))&= r(g(t^5))+r(t^3g(t^5))\\
&=r(g(t^5))+r(t^{15}g(t^5))\\
&\geq r(g(t^5))+r(t^3g(t))=8+2=10.
\end{array}
$$
\item[$n=1$.] This last part cannot be split into more
managable pieces using automorphisms that multiply $t$ by
a root of unity. However, it is possible to use the fact
that in the surface corresponding to $y^2=x^3+tg(t^{10})$,
the fibres over $0$ and over $\infty$ are the same.
To obtain somewhat simpler formulas, note that after scaling
$x,y,t$ by some non-zero elements of $\Qaf$ one can
assume $g(s)=s+1$; we will do this here. An alternative
equation for the one above (replace $x,y$ by $x/t^2, y/t^3$,
respectively) is now $y^2=x^3+t^5+1/t^5$. This is clearly
invariant under $t\mapsto 1/t$. Hence with $u=t+1/t$,
we have a quadratic extension $\Qaf(t)/\Qaf(u)$ and
a curve $E/\Qaf(u)$, for which we are interested in
the points over the quadratic extension. As before,
these points can be split (upto groups of finite index)
into $E(\Qaf(u))$ and $E'(\Qaf(u))$. Here $E'$ is the
corresponding quadratic twist.

Since $t^5+1/t^5=u^5-5u^3+5u$, one finds that $E(\Qaf(u))$
corresponds to the sections of a rational elliptic surface
without reducible fibres. In particular it has rank $8$.
The quadratic extension $\Qaf(t)/\Qaf(u)$ is obtained
by adjoining a square root of $u^2-4$. Hence one concludes
$$
\mbox{ }\hspace{15pt}
r(tg(t^{10}))=r(t^5-5t^3+5t)+r((t^2-4)^3(t^5-5t^3+5t))\geq 8.$$ 
\end{enumerate} 
\end{enumerate}
Summing all lower bounds described above, one finds 
$60$ independent points on Shioda's example, and all these
points can be interpreted in terms of rational elliptic
surfaces.

It remains to prove the assertion about the remaining
rank $68-60=8$ part. One can of course use Shioda's
algorithm from \cite[Thm.~1]{Shi86} to find all other
contributions $r(t^ag(t^b))$. This results in
the equality $r(tg(t^{10}))=16$. Since we only used
the lower bound $8$, this is where the remaining
sections may be found.

Alternatively, this is seen as follows. Start with the
$K3$-surface $Y$ corresponding to $y^2=x^3+s^5g(s^2)$.
Since the elliptic fibration coming from this equation
has $2$ fibres of type $II^*$, the N\'{e}ron-Severi rank
$\rho(Y)$ of $Y$ is at least $2+2\cdot 8=18$.
Substituting $s=t^6,x=t^{8}\xi,y=t^{12}\eta$ yields
a finite rational map $X\rightarrow Y$, with $X$ the
$K3$-surface corresponding to $\eta^2=\xi^3+tg(t^{10})$.
In particular, using Inose's result from \cite{In} once
more, one concludes $\rho(X)=\rho(Y)\geq 18$.
Since the obvious elliptic fibration on $X$ contains
no reducible fibres, this implies $r(tg(t^{10}))\geq 16$.
Using 
$$\begin{array}{rl}
r(tg(t^{10}))&=r(t^5-5t^3+5t)+r((t^2-4)^3(t^5-5t^3+5t))\\
&= 8+r((t^2-4)^3(t^5-5t^3+5t))
\end{array}$$
as was shown above,
one finds that the remaining $8$ sections are in fact
on the elliptic $K3$-surface with equation
$(t^2-4)y^2=x^3+t^5-5t^3+5t$. This finishes the proof
of proposition~\ref{rk36}.
\hfill{$\Box$}
\par
\noindent{\sl Remark $1$.} It may be interesting to note that in
\cite{Fa}, Fastenberg proved that under some conditions which
do not hold in our examples, an elliptic curve $E/\C(t)$
has finite rank over the union of all extensions
$\C(t^{1/n})$. Shioda's upper bound $68$ given in
proposition~\ref{rk36} shows that this is
also true for the curve given by $y^2=x^3+at+b$, provided that
$ab\neq 0$.
\par
\noindent{\sl Remark $2$.} Many variations on these ideas are
possible. For instance, starting from a quadratic polynomial
$f(t)$ such that $r(f(t^6))=18$ (theorem~\ref{main} shows an
example), one can use Inose's result to conclude that
$r(tf(t^5))=18$ as well, and $r(t^3f(t^2))=r(t^5f(t^2))=6$.
Combining as before, one finds $r(f(t^{60}))\geq 54$, where
all these points come from rational or from $K3$-surfaces.

\subsection{Ranks over the rational numbers}
If one restricts to curves $E/\Q(t)$, lemma~\ref{qrank} implies
that the possible rank here is at most half the rank
over $\Qaf(t)$. However, it seems quite hard to get
anywhere near such a bound. Mestre \cite{Me} published an
example of a $j=0$-curve which cannot be defined over $\Q$
and which has $7$ independent points over $\Q(t)$. His
curve is not of the form $y^2=x^3+f(t^6)$ as is studied here.
We tried to find some high rank examples of this form,
by demanding that various associated eigenspaces contain
$\Q(t)$-rational points. The description by Bremner \cite{Br}
of certain classes of curves with positive $\Q(t)$-rank
is useful here. Nevertheless, we had very limited success.

The curve given by $y^2=x^3+t^{12}-26t^6-343$ contains $4$
independent points over $\Q(t)$, namely with
$x$-coordinates $8$, $(-t^6+49)/t^2$, $t^6+7$ and
$(9t^{12}+86t^6+49)/(16t^2)$. This curve is
birational to the one with equation $F(X,Y)=1$, for
$$F(X,Y)=Y^3+(X-Y)(X-2Y)(X-(1-t^3)Y/2).$$
In terms of the latter equation, we consider rational
points with $Y=0,1$ or with $X=0$.

The curve with equation $y^2=x^3+t^{18}+2973t^{12}+369249t^6+11764900$
contains $5$ independent $\Q(t)$-rational points. This example
already appeared in \cite{S-T}, where it was explained how
it is obtained by the same construction that Mestre used when
he finds his rank $\geq 7$ example.

\end{document}